\documentclass[11pt,a4paper,twoside,final]{scrartcl}
\usepackage{a4wide}
\usepackage{amsfonts}
\usepackage{amsmath}
\usepackage{amssymb}
\usepackage{amsthm}
\usepackage{mathtools}
\usepackage[mathscr]{euscript}
\usepackage[numbers]{natbib}
\usepackage{dsfont}
\usepackage{booktabs}
\usepackage{algorithm}
\usepackage{algpseudocode}
\usepackage{pgfplots}
\usepackage{pgfplotstable}
\usepackage{graphicx}
\usepackage{soul}
\usepackage{color}
\usepackage{colortbl}
\usepackage{multirow}
\usepackage{boxedminipage}
\usepackage{enumerate}
\usepackage{collcell}
\usepackage{caption}
\usepackage{subcaption}
\usepackage{rotating}
\usepackage{tablefootnote}
\usepackage{hyperref}
\usepackage{csquotes}
\usepackage{threeparttable}

\usepackage{todonotes}

\newcommand{\A}{\ensuremath{\mathcal{A}}}
\newcommand{\C}{\ensuremath{\mathbb{C}}}

\newcommand{\N}{\ensuremath{\mathbb{N}}}

\newcommand{\R}{\ensuremath{\mathbb{R}}}

\renewcommand{\b}{\boldsymbol}

\newcommand{\ghat}{{\ensuremath{\hat{g}}}}

\newcommand{\abs}[1]{\left\vert #1\right\vert}
\newcommand{\norm}[1]{\left\Vert #1\right\Vert}

\newcommand{\diff}{\ensuremath{\mathrm{d}}}

\DeclareMathOperator*{\argmin}{arg\,min}

\definecolor{color1}{HTML}{E69F00}
\definecolor{color2}{HTML}{56B4E9}
\definecolor{color3}{HTML}{009E73}
\definecolor{color4}{HTML}{F0E442}
\definecolor{color5}{HTML}{0072B2}
\definecolor{color6}{HTML}{D55E00}
\definecolor{color7}{HTML}{CC79A7}

\newtheorem{theorem}{Theorem}[section]
\newtheorem{lemma}[theorem]{Lemma}
\newtheorem{remark}[theorem]{Remark}
\newtheorem{generalisation}[theorem]{Generalisation}
\newtheorem{definition}[theorem]{Definition}
\newtheorem{example}[theorem]{Example}
\newtheorem{corollary}[theorem]{Corollary}
\newtheorem{proposition}[theorem]{Proposition}

\usepgfplotslibrary{statistics}
\usetikzlibrary{patterns}
\pgfplotsset{compat=newest}
\usetikzlibrary{calc}
\usetikzlibrary{positioning,arrows.meta}

\pgfplotsset{
        my boxplot style/.style={
            boxplot prepared,
            draw=black,
            solid,
            fill=white,
            mark=*,
            every mark/.append style={
                fill=gray,
            },
        },
    }

\pgfplotsset{
    boxplot prepared from table/.code={
        \def\tikz@plot@handler{\pgfplotsplothandlerboxplotprepared}%
        \pgfplotsset{
            /pgfplots/boxplot prepared from table/.cd,
            #1,
        }
    },
    /pgfplots/boxplot prepared from table/.cd,
        table/.code={\pgfplotstablecopy{#1}\to\boxplot@datatable},
        row/.initial=0,
        make style readable from table/.style={
            #1/.code={
                \pgfplotstablegetelem{\pgfkeysvalueof{/pgfplots/boxplot prepared from table/row}}{##1}\of\boxplot@datatable
                \pgfplotsset{boxplot/#1/.expand once={\pgfplotsretval}}
            }
        },
        make style readable from table=lower whisker,
        make style readable from table=upper whisker,
        make style readable from table=lower quartile,
        make style readable from table=upper quartile,
        make style readable from table=median,
        make style readable from table=lower notch,
        make style readable from table=upper notch,
        make style readable from table=draw position,
        make style readable from table=box extend,
}
\makeatother

\numberwithin{equation}{section}
\numberwithin{table}{section}
\numberwithin{figure}{section}

\mathtoolsset{showonlyrefs=true}  

\title{Learning solution operators of PDEs with sparse approximation methods}
\author{Sebastian Neumayer \and Daniel Potts \and Fabian Taubert}
\date{\today}

\begin{document}

\maketitle

\begin{abstract}
We investigate the approximation of solution operators for partial differential equations (PDEs) using sparse high-dimensional techniques.
Building on a dimension-incremental framework, we combine product basis expansions with sparse recovery methods, specifically orthogonal matching pursuit (OMP), to substantially reduce the required sample size compared with a previously considered cubature-based approach.
We evaluate the resulting method numerically on several examples, comparing it against both cubature-based sparse approximation and Fourier neural operators in terms of accuracy, runtime, and sample size.
The experiments show that our approach considerably reduces the number of required PDE solves relative to its predecessor while maintaining competitive accuracy, particularly when the solution admits a sparse representation in the chosen basis. 
Furthermore, the recovered sparse index sets yield interpretable insights into the relevant variables and parameter interactions.

\medskip

\end{abstract}

\section{Introduction}
Parameter-dependent partial differential equations (PDEs) arise naturally across a broad spectrum of scientific and engineering applications, including forward simulations, control theory, uncertainty quantification, and inverse problems.
There, the repeated numerical solution of the underlying PDEs often represents the dominant computational cost, particularly for high-dimensional parameter spaces.
Consequently, the efficient approximation of the solution operator, which maps parameter inputs such as initial values or forcing terms  to the PDE solution, has emerged as a central problem in numerical analysis and scientific computing.

Recent studies indicate that learned approaches, see \cite{BoTo24,KoLaStu24,SuTe26} for an overview, can outperform classical techniques \cite{azizzadenesheli2024neural}.
As the overall computational effort is often dominated by the number of required PDE solves for generating the training samples rather than by the training itself, methods should be sample efficient.
To formalize operator learning in the PDE context, let $\mathcal{F}$ denote a (potentially infinite-dimensional) parameter space, $\mathcal{U}$ a target space of functions on some domain $\Omega \subset \mathbb{R}^d$, and $G \colon \mathcal{F} \to \mathcal{U}$ the unknown solution operator.
As in reduced-order modeling \cite{BhHoKoSt21} and consistent with the finite capacity of the hidden representations used by neural operators \cite{kovachki2023neural}, we project the input parameters onto a finite set of basis coefficients $\b a \in \mathbb{R}^n$, yielding a finite-dimensional operator $G_n \colon \mathbb{R}^n \to \mathcal{U}$.
Given training data $\{(\b a_i, \b x_i, u_{\b a_i}(\b x_i))\}_{i=1}^M \subset \mathbb{R}^n \times \Omega \times \mathbb{R}$, consisting of parameter vectors paired with solution evaluations (often with multiple spatial evaluations per parameter $\b a$), the goal is to learn an accurate approximation of $G_n$.
A central difficulty is the curse of dimensionality: the sample and computational complexity may scale exponentially in the input dimension $n$.

To construct a tractable approximation of $G_n$, one must first choose a suitable parametric architecture.
In approximation theory, a standard approach for handling the
infinite-dimensional range space $\mathcal{U}$ is to choose a basis $\{\Phi_j\}_{j \in J} \subset \mathcal{U}$, such as the Fourier or Chebyshev basis, and to work with the truncated expansion
\begin{equation}\label{eq:OperatorApprox}
    G_n[\b a](\b x) \approx \sum_{j \in J} b_j(\b a) \Phi_j(\b x)
    \qquad \forall\, \b x \in \Omega.
\end{equation}
This expansion is also used by the DeepONet approach \cite{LuJiPaZhaKa21}.
However, we restrict \eqref{eq:OperatorApprox} to a fixed, non-learnable basis.
Then, parameterizing $b \colon \mathbb{R}^n \to \mathbb{R}^{|J|}$ by a neural network $b(\cdot\,;\b c)$ with parameters $\b c \in \mathbb{R}^P$, as in DeepONet's branch network, leads to the training problem
\begin{equation}\label{eq:NN}
    \argmin_{\b c \in \R^{P}} \sum_{i=1}^M \biggl\vert
   \sum_{j \in J} b_j(\b a_i; \b c) \Phi_j(\b x) - u_{\b a_i}(\b x_i) \biggr\vert^2.
\end{equation}
By expanding a neural operator \cite{kovachki2023neural} output in the basis $\{\Phi_j\}_{j \in J}$, we can formally interpret its learning within the same framework.
While expressive, such nonlinear architectures typically require many training samples and offer limited theoretical guarantees.
In contrast, we adopt a simple linear feature model of the form $b_j(\b a; \b c) = \sum_{k \in K} c_{k,j} \Psi_k(\b a)$, where $\{\Psi_k\}_{k \in K}$ is a fixed dictionary of features $\Psi_k \colon \mathbb{R}^n \to \mathbb{R}$ and the coefficients $c_{k,j} \in \mathbb{R}$ are the learnable parameters.
Setting $\b k = (k,j) \in \Gamma \coloneqq K \times J$ and $T_{\b k}(\b a, \b x) \coloneqq \Psi_k(\b a) \Phi_j(\b x)$, the resulting approximation \eqref{eq:OperatorApprox} can be compactly written as $G_n[\b a](\b x) \approx \sum_{\b k \in \Gamma} c_{\b k} T_{\b k}(\b a, \b x)$.
Then, the generic operator learning task \eqref{eq:NN} reduces to the least-squares problem
\begin{equation}\label{eq:LeastSquares}
    \argmin_{\b c \in \mathbb{R}^{\vert \Gamma \vert}} \sum_{i=1}^M \biggl\vert \sum_{\b k \in \Gamma} c_{\b k} T_{\b k}(\b a_i, \b x_i) - u_{\b a_i}(\b x_i) \biggr\vert^2,
\end{equation}
which is theoretically well analyzed. 
However, since the number of unknowns $|\Gamma|$ typically far exceeds the number of training samples $M$, solving \eqref{eq:LeastSquares} directly leads to severe overfitting. 
Sparse approximation methods \cite{KuRa08,ChCoSchw14,DuTeUl16,KaUlVo19} address this by assuming that $G_n$ can be represented accurately using only a small number of nonzero coefficients $c_{\b k}$.
This assumption is not merely a regularization device: the analysis of specific PDEs shows that the analytical solutions indeed satisfy such sparsity constraints \cite{PoTa25}.
Together, these observations motivate the constrained training problem
\begin{equation}\label{eq:SparseRec}
    \argmin_{I \subset \Gamma:\, |I| \leq s}\; \min_{\b c \in \mathbb{R}^{|I|}} \sum_{i=1}^M \biggl\vert \sum_{\b k \in I} c_{\b k} T_{\b k}(\b a_i, \b x_i) - u_{\b a_i}(\b x_i) \biggr\vert^2.
\end{equation}
The problem \eqref{eq:SparseRec} can be solved, e.g., with greedy algorithms such as orthogonal matching pursuit (OMP) \cite{PaReKr93}.
Under suitable conditions on the feature dictionary and the sample distribution, robust recovery guarantees for \eqref{eq:SparseRec} are available; see, e.g., \cite{Zha11,FoRa13,MoNePoSoUl25}.

While sparse recovery addresses the overfitting problem, the cardinality of the search space $\Gamma$ can itself be prohibitively large.
For a product-type basis for $\Phi_j$ and $\Psi_k$ with, e.g., $64$ basis functions per direction, the number of candidate indices already grows to $|\Gamma| \approx 64^{n+d}$, rendering an exhaustive search computationally infeasible.
As a remedy, the dimension-incremental framework of \cite{PoVo14,KaPoVo17,KaPoTa23} constructs the support set $I$ dimension by dimension: beginning with a one-dimensional candidate set, it successively extends the current support rather than searching $\Gamma$ directly.
The extension step involves solving several sparse recovery problems of the form \eqref{eq:SparseRec}, each with a tractable number of candidates.
In this way, the method can efficiently exploit the underlying sparsity structure while considerably reducing both computational and sample complexity compared to an exhaustive search.
A recent alternative for managing large (continuous) candidate sets is generative feature training \cite{HerNeu2025}, which adaptively resamples the candidates for \eqref{eq:SparseRec} in place of a greedy approach.

\paragraph{Contributions}
Our contributions can be summarized as follows.

\begin{itemize}
    \item We combine the dimension-incremental support detection framework  of \cite{PoVo14,KaPoVo17} with OMP-based sparse recovery \cite{MoNePoSoUl25} to compute sparse approximations of PDE solution operators from carefully selected samples.
    Compared to cubature-based recovery \cite{KaPoTa23}, this significantly reduces the sample size and consequently the number of PDE solves.
    
    \item By design, the cubature-based recovery in \cite{PoTa25} requires an independent PDE solve for each sample.
    In contrast, sparse recovery methods can naturally incorporate multiple samples with different spatial locations obtained from a single solve, which improves the computational efficiency even further. 
    
    \item We numerically validate the proposed methods for several PDEs, comparing approximation quality, runtime, and sampling complexity against cubature-based sparse approximation methods and Tensorized Fourier Neural Operators (TFNOs) as a baseline.
    
    \item Unlike purely data-driven black-box approaches, our approach provides interpretable information about the parameter dependence of the solution via the detected index sets and their interaction structure.
\end{itemize}

\begin{figure}[p]
\centering
\begin{tikzpicture}[
    box/.style={
        draw,
        align=center,
        rounded corners=1pt,
        inner sep=5pt,
        minimum height=1.0cm
    },
    arr/.style={-{Latex[length=3mm,width=2.2mm]}, thick},
    node distance=1.2cm and 1.8cm
]

\node[box] (Gn) {
    Finite-dimensional operator\\
    $G_n\colon \mathbb R^n \to \mathcal U$
};

\node[box, below left=0.5cm and -1.0cm of Gn] (NNform) {
    Formal neural operator expansion\\[.5ex]
    $\displaystyle
    G_n[\b a](\b x)
    \approx
    \sum_{j\in J} b_j(\b a; \b c)\Phi_j(\b x)$
};

\node[box, below right=0.5cm and -1.0cm of Gn] (Sparseform) {
    $\displaystyle
    G_n[\b a](\b x)
    \approx
    \sum_{\b k\in I} c_{\b k}T_{\b k}(\b a,\b x)$
};

\node[below=4.5cm of Gn] (ways) {};

\node[box, left=4cm of ways] (tfno) {
    train a \texttt{TFNO} \\
    baseline \cite{kossaifi2024multigrid}\\
    by solving \eqref{eq:NN}
};

\node[box, left=0cm of ways] (ompp) {
    solve \eqref{eq:SparseRec}\\
    via \texttt{OMP+}\\
    (dim.-incremental)
};

\node[box, right=0cm of ways] (omp) {
    solve \eqref{eq:SparseRec}\\
    via \texttt{OMP}\\
    (dim.-incremental)
};

\node[box, right=4cm of ways] (cub) {
    compute $c_{\b k}$\\
    by cubature rule\\
    for comparison \cite{PoTa25}\\
    (dim.-incremental)
};

\node[box, below=11cm of Gn] (data) {
    solve $M$ PDEs to obtain training data\\
    $(\b a_i,\b x_i,u_{\b a_i}(\b x_i))$, $i=1,\ldots,M$
};

\node[box, above left=1.5cm and -3cm of data] (virtual) {
    generate additional data for fixed $\b a_i$\\
    (freely available from PDE solver)\\[2pt]
    $(\b a_i,\b x_{i,j},u_{\b a_i}(\b x_{i,j}))$
};

\node[circle, fill=black, inner sep=1.4pt, below=0.75cm of Sparseform](j1) {};
\node[circle, fill=black, inner sep=1.4pt, above=0.75cm of data](j2) {};
\node[circle, fill=black, inner sep=1.4pt, above=0.5cm of virtual](j3) {};

\coordinate(marker_tfno) at (NNform |- j1);

\draw[arr] (Gn) -| (NNform);
\draw[arr] (Gn) -| (Sparseform);

\draw[thick] (NNform) -- (marker_tfno);
\draw[arr] (marker_tfno) -| (tfno);

\draw[thick] (Sparseform) -- (j1);
\draw[arr] (j1) -| (ompp);
\draw[arr] (j1) -| (omp);
\draw[arr] (j1) -| (cub);

\draw[thick] (data) -- (j2);
\draw[arr] (j2) -| (virtual);
\draw[arr] (j2) -| (omp);
\draw[arr] (j2) -| (cub);

\draw[thick] (virtual) -- (j3);
\draw[arr] (j3) -| (tfno);
\draw[arr] (j3) -| (ompp);

\end{tikzpicture}
\caption{Schematic comparison of the considered approximation strategies.}
\label{fig:method_overview}
\end{figure}
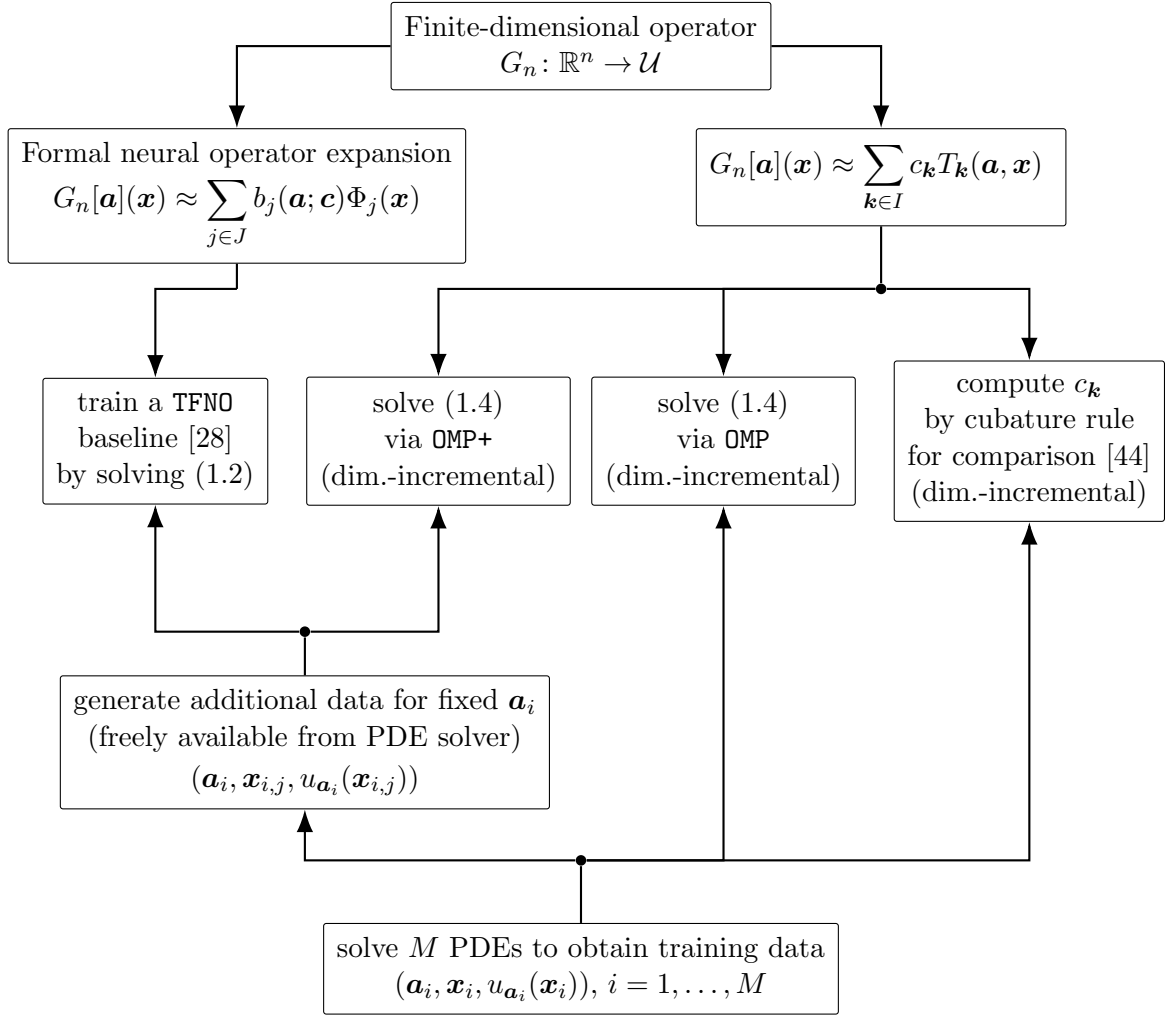

Figure \ref{fig:method_overview} provides a schematic overview of the proposed approximation strategies together with the methods used for comparison.
Moreover, Table \ref{tab:overview} summarizes typical parameter sizes occurring in the numerical experiments of Section \ref{sec:numerics}, as well as the corresponding numbers of PDE solves required by the different approaches.
The Python implementation and numerical experiment scripts are publicly available at \cite{gitNABOPB}.

\begin{table}[p]
\centering
\begin{threeparttable}
\begin{tabular}{lccc}
\toprule
quantity & Heat eq. & Burgers' eq. & Parametric diffusion eq.\\
\midrule
overall dimension $d_\mathrm{ex}$ & $11$ & $10$ & $22$ \\
amount of candidates $|\Gamma|$ & $8.75 \cdot 10^{19}$ & $1.35 \cdot 10^{18}$
    & $7.55\cdot 10^{17}$\tnote{$*$} \\
target sparsity $s$ & $1000$ & $1000$ & $1000$ \\
\texttt{TFNO} parameters & $529\,233$ & $81\,049$ & $529\,233$ \\
\hline
PDE solves & & & \\
\texttt{OMP}   & $285\,000$ & $255\,000$   & $615\,000$ \\
\texttt{OMP+}  & $35\,006$  & $35\,002$    & $79\,006$  \\
\texttt{cMR1L} & $702\,054$ & $5\,706\,215$\tnote{$\dagger$} & $3\,683\,946$ \\
\texttt{TFNO}  & $40\,500$  & $40\,500$    & $86\,000$  \\
\bottomrule
\end{tabular}
\begin{tablenotes}
  \item[$*$] Additional assumption: $\|\b k\|_0 \leq 7 \quad \forall \b k \in \Gamma$.
  \item[$\dagger$] Here, $s=500$ and a smaller  $\Gamma$ were used due to computational restrictions.
\end{tablenotes}
\end{threeparttable}
\caption{Parameter counts and number of PDE solves for our
  experiments in Section~\ref{sec:numerics}.}
\label{tab:overview}
\end{table}

\paragraph{Outline}

The remainder of this work is organized as follows.
Section~\ref{sec:theory} discusses the parametric representation of PDE solution operators and their sparse approximation.
Subsequently, we recall the dimension-incremental support detection framework and combine it with OMP-based coefficient recovery and an improved sampling strategy for PDE-based problems.
Section~\ref{sec:numerics} presents numerical experiments for the heat equation, Burgers' equation, and a parametric diffusion equation, comparing approximation quality, runtime, and sampling complexity against cubature-based approaches and TFNOs.
Finally, Section~\ref{subsec:discussion} discusses the numerical results and the trade-offs between approximation quality and computational cost.

\paragraph{Related work}
The approximation of (infinite-dimensional) operators from training data has recently attracted significant attention, with approaches broadly falling into three categories:
linear expansions, neural operators, and (random) kernel methods.

From the viewpoint of approximation theory, operator learning can be cast as finding the coefficients of a basis expansion \cite{deHoop2023convergence, DuTeUl16}.
This perspective is closely related to sparse polynomial approximation \cite{ChCoSchw14,RaWa12}: if the solution operator admits a sparse representation in a suitable basis, only a small number of coefficients need to be recovered \cite{BacCohDa18}.
To this end, a PDE-oriented recovery strategy was investigated in \cite{PoTa25}.
Alternatively, one can directly learn mappings between the basis expansion coefficients of input and output functions using
neural networks \cite{kissas2022nomad,ingebrand2025basis} or approximate the Green's function instead \cite{boulle2022greens}.

Neural operators are typically trained in a supervised setting using (potentially physics-informed \cite{li2024physicsinformed}) loss functions.
Prominent architectures include DeepONets \cite{LuJiPaZhaKa21,wang2022improved} and Fourier Neural Operators (FNOs) \cite{LiEtAl21}, many of which are implemented as part of the NeuralOperator library \cite{gitNeuralOperator}.
Spectral neural operators and related architectures replace purely Fourier-based representations with polynomial or spectral basis expansions \cite{FaOs24,SchiaEtAl24}, a perspective closely related to this work.
Other recent directions include the incorporation of attention mechanisms \cite{hao2023gnot,li2023transformer} and the aim for foundation models \cite{sun2025towards}.
Comprehensive surveys covering practical and theoretical aspects of neural operators can be found in \cite{BoTo24,KoLaStu24,DuKoAn26}.

Operator-valued kernels provide a principled framework \cite{KaDuPrCaRaAu16}, though the direct implementation is computationally prohibitive in general.
Random feature approaches \cite{nelsen2024operator} offer a scalable alternative: by fixing the random parameters during training, the optimization reduces to a linear problem.
Interestingly, kernel-based methods can achieve accuracy competitive with neural operators \cite{BaDaHoOw24}.
A related Gaussian process framework was recently proposed in \cite{mora2025operator}.

Finally, we remark that there is also growing interest in topics such as error estimates \cite{lanthaler2022error,kovachki2023neural}, discretization consistency \cite{bartolucci2023representation} and sample complexity \cite{kovachki2024data} of neural operators.

\section{Theory}
\label{sec:theory}

The presentation of the general framework follows~\cite{PoTa25}, where operator learning for partial differential equations (PDEs) is investigated from a sparse high-dimensional approximation perspective.
For the purposes of this paper, we adapt their approach and incorporate dedicated sparse recovery algorithms.
While we restrict ourselves to time-dependent PDEs as a running example below, we stress that the proposed method applies to a broad class of operators, which also includes stationary and parametric PDEs.
\subsection{Time-dependent partial differential equations}
\label{subsec:problem_setting}
Let $\Omega \subset \mathbb{R}^d$ be a bounded spatial domain and $T>0$.
Further, assume that we are given function spaces $\mathcal{U}$ and $\mathcal{F}$ defined on the space-time domain $\Omega \times (0,T)$, a function space $\mathcal H$ on $\partial \Omega \times (0,T)$, and a function space $\mathcal U_0$ on $\Omega$.
Here, the pairs $(\mathcal U, \mathcal U_0)$ and $(\mathcal U, \mathcal H)$ should be compatible in the sense that $u(\cdot,0) \in \mathcal U_0$ and $u|_{\partial\Omega \times (0,T)} \in \mathcal H$ for any $u \in \mathcal U$.
Within this setting, we consider differential equations of the form
\begin{align}\label{eq:general}
L u = f \qquad \qquad \text{s.t. } u|_{\partial\Omega \times (0,T)} = h \text{ and } u(\cdot,0)  = u_0,
\end{align}
where $L \colon \mathcal{U}\to\mathcal{F}$ is a (possibly nonlinear) differential operator, $f\in\mathcal{F}$ is an inhomogeneity, and $u_0 \in \mathcal U_0$ is an initial value.
Throughout, we assume that \eqref{eq:general} is well-posed in the sense that for every $(f,h,u_0) \in \mathcal F \times \mathcal{H} \times \mathcal{U}_0$, there exists a unique solution $u \in \mathcal{U}$ depending continuously on the data $(f,h,u_0)$.
Then, \eqref{eq:general} induces a solution operator
\begin{equation}\label{eq:solutionOp}
G\colon \mathcal{F} \times \mathcal{H} \times \mathcal{U}_0 \to \mathcal{U},
\qquad (f,h,u_0) \mapsto u.
\end{equation}
In practice, often only a subset of the data $(f,h,u_0)$ is varied, while the remaining components are fixed.
In particular, a common situation is that the solution operator $G$ reduces to a mapping from the initial data $u_0$ to the solution $u$, with $f$ and $h$ being fixed.
Then, our goal is to approximate $G$ by a surrogate model that enables fast evaluations for new input data $u_0$ without the need for solving the PDE \eqref{eq:general}.

\begin{remark}
\label{rem:parameters}
For the Poisson equation as discussed in \cite[Sec.~3]{PoTa25}, the solution would depend only on the spatial variable $\b x$ and the inhomogeneity $f$.
Such variations mainly influence the dimensionality of the product spaces in \eqref{eq:solutionOp} and the choice of discretization, while leaving the approximation procedure for $G$ described below unchanged.
\end{remark}

\subsection{Parametric operator representation and sparse approximation}
\label{subsec:sparse_approximation}
In order to obtain a finite-dimensional representation of the domain of $G$, we expand the initial condition $u_0$ as
\begin{align}
u_0(\b x) \approx \sum_{j=1}^n a_j A_j(\b x),
\qquad \b x \in \Omega,
\label{eq:u0_approx}
\end{align}
with basis functions $A_j \colon \Omega \to \R$, $j=1,\ldots,n$, and coefficients $\b a = (a_1,\dots,a_n)\in[-1,1]^n$.
Note that the chosen functions $A_j$ should lead to a sufficiently accurate approximation in $\mathcal U_0$.
Based on \eqref{eq:u0_approx}, we get the induced solution operator $G_n \colon \mathbb [-1,1]^n \to \mathcal U$ with $\b a \mapsto u$.

In general, any appropriate discretization $G_n$ of $G$ can be identified with a function defined on an extended bounded product domain \smash{$\mathcal{D} \coloneqq \prod_{j=1}^{d_\mathrm{ex}} \mathcal{D}_{j} \subset \mathbb{R}^{d_\mathrm{ex}}$}.
In the case of \eqref{eq:u0_approx}, we get $\mathcal D = \Omega \times [0,T] \times [-1,1]^n$.
For the next step, we assume that the induced function admits an expansion in a bounded orthonormal product basis (BOPB) $\{\Phi_{\b k} : \b k \in \N^{d_\mathrm{ex}}\} \subset L_2(\mathcal{D},\mu)$, where $\mu$ is a product measure.
Recall that a BOPB consists of tensor-products $\Phi_{\b k}(\b y) = \prod_{j=1}^{d_\mathrm{ex}} \phi_{j,k_j}(y_j)$.
More details can be found, e.g., in \cite{KaPoTa23}.

For the operator $G_n$ induced by \eqref{eq:u0_approx}, the BOPB expansion reads as
\begin{align}
u(\b x, \tau; \b a)
=
\sum_{\b k \in \N^{d_\mathrm{ex}}} c_{\b k}\,\Phi_{\b k}(\b x, \tau, \b a),
\qquad
(\b x, \tau, \b a)\in \mathcal D,
\label{eq:basexp_of_u}
\end{align}
with $c_{\b k}\in \C$ being the corresponding basis coefficients
\begin{align}\label{eq:integral_ck}
\hat c_{\b k} = \int_\mathcal{D} g(\b x, \tau, \b a) \overline{\Phi_{\b k}(\b x, \tau, \b a)}\mathrm d (\b x, \tau, \b a).
\end{align}
Then, our aim is to accurately approximate $u$ (if possible) by a sparse truncation
\begin{align}
S_I u(\b x, \tau; \b a)
=
\sum_{\b k\in I} c_{\b k}\,\Phi_{\b k}(\b x, \tau, \b a),
\label{eq:basexp_of_u_trunc1}
\end{align}
where $I \subset \N_0^{d_\mathrm{ex}}$ is a finite index set.
Since the exact $c_{\b k}$ are unknown, we replace them by approximations $\hat u_{\b k}$ (for which we outline a computation procedure below) and obtain
\begin{align}
S_I^\A u(\b x, \tau, \b a)
=
\sum_{\b k\in I} \hat u_{\b k}\,\Phi_{\b k}(\b x, \tau, \b a).
\label{eq:basexp_of_u_trunc}
\end{align}
In summary, \eqref{eq:basexp_of_u_trunc} yields a continuous approximation of $G$ by first projecting $u_0$ onto its finite coefficient vector $\b a$ and then evaluating $S_I^\A u(\cdot,\cdot,\b a)$.
The resulting approximation error can be split into a truncation part and a coefficient recovery part.
More precisely, using the uniform boundedness of the basis functions $\Phi_{\b k}$, we obtain
\begin{align}
\norm{u-S_I^\A u}_\infty
&\leq
\norm{u-S_Iu}_\infty + \norm{S_Iu-S_I^\A u}_\infty \notag\\
&\leq
\sum_{\b k\notin I} \abs{c_{\b k}} \norm{\Phi_{\b k}}_\infty
+
\sum_{\b k\in I} \abs{c_{\b k}-\hat u_{\b k}} \norm{\Phi_{\b k}}_\infty \notag\\
&\leq
\sup_{\b k \in I} \norm{\Phi_{\b k}}_\infty \biggl(
\sum_{\b k\notin I} \abs{c_{\b k}}
+
\sum_{\b k\in I} \abs{c_{\b k}-\hat u_{\b k}}
\biggr).
\end{align}
Hence, a good approximation requires both identifying a suitable index set $I$ and an accurate method for estimating the coefficients $c_{\b k}$ from samples $u(\b x_j,\tau_j,\b a_j)$.
For our running example, we keep in mind that samples are obtained by numerically solving the PDE \eqref{eq:general} for the parameters $\b a_j$, which makes sampling costly.

Given a finite index set $I$, the corresponding coefficients $c_{\b k},\, \b k \in I,$ can be estimated using a cubature formula based on \eqref{eq:integral_ck}.
Alternatively, they may be obtained by solving the linear system that arises from inserting the samples $u(\b x_j,\tau_j,\b a_j)$ into \eqref{eq:basexp_of_u}, for instance via a standard least squares approach.
However, since the index set $I$ is not known a priori, one typically has to work with the product space $\Gamma \supset I$ instead.
Without imposing strong a priori assumptions, such search spaces grow rapidly in high-dimensional settings and render the aforementioned approaches computationally infeasible.
We therefore turn to the problem of determining a suitable index set $I \subset \Gamma \subset \N_0^{d_\mathrm{ex}}$.

\begin{remark}
\label{rem:transform}
In our implementation, we transform the problem-specific product domain $\mathcal{D}$ to the hypercube $[-1,1]^{d_\mathrm{ex}}$ via problem-dependent affine maps, specified for each PDE individually.
After this transformation, Chebyshev polynomials can serve as a BOPB for the function induced by the operator $G_n$.
\end{remark}

\subsection{Dimension-incremental support detection and sparse recovery}
\label{subsec:dim_incremental}
Since the candidate set $\Gamma$ has product structure, its cardinality scales as $10^{d_\mathrm{ex}}$ even with only 10 basis functions per dimension.
For $d_\mathrm{ex} = 12$, this already yields $10^{12}$ candidate
index tuples, making an exhaustive search for a suitable $I$ in general computationally infeasible.

As an alternative, we employ the dimension-incremental support detection strategy proposed in \cite{KaPoTa23}.
The main idea is to build $I$ progressively, starting from one-dimensional projections and then combining the detected candidates across dimensions.
Below we provide a brief summary; see Algorithm~\ref{alg:main} and \cite{KaPoTa23} for further details and a detection guarantee for a suitable index set $I$.
To simplify the notation, we switch from the problem-specific function $u$ to a generic function $g$, and correspondingly collect all variables into a single $d_\mathrm{ex}$-dimensional input vector.

Let $g\colon \R^{d_\mathrm{ex}}\to\C$ be an (unknown) target function with ($s$-sparse) expansion
\begin{equation}\label{eq:expansion}
g = \sum_{\b k\in I} \hat g_{\b k}\Phi_{\b k},
\end{equation}
where $I$ with $\vert I \vert =s$ is contained in a prescribed finite candidate space $\Gamma \subset \N_0^{d_\mathrm{ex}}$.
Typically, $\Gamma $ is simply a cube in $\N_0^{d_\mathrm{ex}}$.
For each coordinate $t\in\{1,\dots,d_\mathrm{ex}\}$, we consider the one-dimensional projected candidate set $\mathcal{P}_{\{t\}}(\Gamma) \coloneqq \{ l \in \N_0 \,|\,\exists \b k \in \Gamma: k_t=l\}$ and for each $l \in \mathcal{P}_{\{t\}}(\Gamma)$ the so-called projected coefficients
\begin{align}
\hat g_{\{t\},l}(\boldsymbol{\theta})
\coloneqq
\int_{\mathcal{D}_t} g_{\{t\}}(\xi,\boldsymbol{\theta})
\overline{\phi_{t,l}(\xi)}
\,\diff \xi = \sum_{\b k\in I : k_t = l} \hat g_{\b k} \prod_{j \neq t} \phi_{j, k_j}(\theta_j),
\label{eq:proj_coef_1d_new}
\end{align}
where $\boldsymbol{\theta} \in \R^{d_\mathrm{ex}-1}$ is an anchor and $g_{\{t\}}(\xi,\boldsymbol{\theta})$ denotes the value of $g$ with $\xi$ in the $t$-th dimension and $\boldsymbol{\theta}$ for the remaining coordinates.
Due to the representation \eqref{eq:expansion}, at most $s$ coefficients $\hat g_{\{t\},l}(\boldsymbol{\theta})$ should be nonzero.
Thus, among all coefficients $l \in \mathcal{P}_{\{t\}}(\Gamma)$, we retain the $s$ largest ones (in terms of $|\hat g_{\{t\},l}(\boldsymbol{\theta})|$) exceeding a prescribed threshold, and discard all others as negligible.
All detected candidates $l$ are collected in the temporary index sets $I_{\{t\}}$.

\begin{remark}\label{rem:detection_iterations}
The projected coefficients $\hat g_{\{t\},l}(\boldsymbol{\theta})$ are an indicator for the relevance of the index $l$. 
However, an unfavorable choice of the anchor $\boldsymbol{\theta}$ may cause cancellations in \eqref{eq:proj_coef_1d_new}, leaving some relevant indices undetected. 
A remedy is to perform $r \in \mathbb{N}$ detection iterations with distinct anchors $\boldsymbol{\theta}$ and to consider the union of the detected indices. 
While lower bounds on $r$ required to achieve a prescribed failure probability can be derived analogously to \cite[Thm.~3.5]{KaPoTa23}, in practice, considerably smaller values of $r$ often suffice for reliable results.

Unfortunately, increasing $r$ raises both the required sample size (and thus the number of PDE solves) and the overall computational effort. 
To mitigate this, the anchors $\boldsymbol{\theta} \in \mathbb{R}^{d_\mathrm{ex}-1}$ can be selected using a distance-maximizing strategy rather than purely random sampling. 
This minimizes the risk of multiple anchors capturing redundant local behavior of $\hat g_{\{t\},l}(\boldsymbol{\theta})$. 
Such local dependencies occur in the heat equation example in Section~\ref{subsec:heat}, where certain solution components decay rapidly over time, rendering anchors with large $t$ blind to them.

\end{remark}

Next, for $t\geq 2$, the candidate sets are combined incrementally.
Given $I_{\{1,\dots,t-1\}}$, one forms the sets
\begin{align}
\mathcal{P}_{\{1,\ldots,t\}}(\Gamma) &\coloneqq \{ \b l \in \N_0^t \,|\,\exists \b k \in \Gamma \subset \N_0^{d_\mathrm{ex}}: \b k_{\{1,\ldots,t\}}= \b l\}\\
K
&\coloneqq
\bigl(I_{\{1,\dots,t-1\}}\times I_{\{t\}}\bigr)
\cap
\mathcal{P}_{\{1,\dots,t\}}(\Gamma), \label{eq:candidate_set}
\end{align}
and investigates the corresponding projected coefficients on $K$ as in the one-dimensional case.
In particular, we again consider the projected coefficients
\begin{align}
\hat g_{\{1,\ldots,t\},\b l}(\boldsymbol{\theta})
\coloneqq\!\!\!
\int\limits_{\prod_{j=1}^t \mathcal{D}_j} g_{\{1,\ldots,t\}}(\b \xi,\boldsymbol{\theta})
\overline{\Phi_{\{1,\ldots,t\},\b l}(\b \xi)}
\,\diff \b \xi = \!\!\!\! \sum_{\b k\in I : \b k
_{\{1,\ldots,t\}}
\!\!= \b l} \hat g_{\b k} \!\!\!\prod_{j \notin \{1,\ldots,t\}} \Phi_{j, k_j}(\theta_j),
\label{eq:proj_coef_hd_new}
\end{align}
estimate their absolute value, retain the $s$ largest ones exceeding the threshold, and collect them in the index set $I_{\{1,\dots,t\}}$.
In this way, higher-order interactions are identified successively until the final index set $I_{\mathrm{est}} \coloneqq I_{\{1,\dots,d_\mathrm{ex}\}}$ is obtained.

In \cite{PoTa25}, the projected coefficients defined in \eqref{eq:proj_coef_1d_new} and \eqref{eq:proj_coef_hd_new} are approximated via cubature rules, which typically require a large number of sampling points $\smash{\b \xi^{(j)}}$ at carefully chosen locations.
Specifically, so-called rank-1 lattices are employed, which provide accurate (and in some cases even exact) cubature for functions with suitable smoothness or frequency structure.
When using large lattices, this enables reliable and efficient recovery of the projected coefficients for all candidate indices, including the ones that are actually negligible.

In contrast to this exhaustive approximation, we propose to directly exploit the system's inherent sparsity.
More precisely, \eqref{eq:expansion} implies that also the projected expansions
\begin{equation}\label{eq:SparseRecProb}
    g_{\{t\}}(\xi,\b \theta) = \sum_{l \in \mathcal{P}_{\{t\}}(\Gamma)} \hat g_{\{t\},l}(\b \theta) \Phi_{\{t\},l}(\xi),
\end{equation}
and similarly the ones with \eqref{eq:proj_coef_hd_new}, are $s$ sparse.
Classical sparse recovery algorithms, such as the orthogonal matching pursuit (OMP) \cite{KuRa08} discussed in Section~\ref{subsec:OMP}, typically require only around $M=5s$ random sample locations $\boldsymbol \xi^{(j)}$, thereby considerably improving the computational efficiency of Algorithm~\ref{alg:main} compared to \cite{PoTa25}.
Since our approach treats $g$ as a black box, any numerical scheme, such as finite elements, finite differences, or spectral methods, can be used to compute the function values $g_{\{1,\ldots,t\}}(\b \xi^{(j)}, \boldsymbol{\theta})$.
In turn, the accuracy of these evaluations directly affects the quality of the computed approximation \eqref{eq:basexp_of_u_trunc}.
Finally, we emphasize that the sampled PDE input parameters $\b a_i$ must conform to the anchor decomposition: each $\b a_i$ is partitioned into active coordinates, sampled randomly, and fixed anchor coordinates.

\renewcommand{\algorithmicrequire}{\textbf{Input:}}
\renewcommand{\algorithmicensure}{\textbf{Output:}}
\begin{algorithm}[t]
    \caption{Dimension-Incremental Index Detection (Simplified)}\label{alg:main}
    \small
    \begin{algorithmic}[1]
        \Require Search space $\Gamma\subset\mathbb{N}^{d_\mathrm{ex}}$, target function $g$, sparsity level $s\in\mathbb{N}$, detection threshold $\delta_+ > 0$, number of detection iterations $r\in\mathbb{N}$.
        
        \Statex \textbf{(Step 1) Single component identification}
        \For{$t=1,\ldots,J$}
            \State Set $I_{\{t\}}\coloneqq \emptyset$.
            \State Draw random sampling nodes \smash{$(\xi^{(j)})_{j=1}^M \subset \mathcal{D}_t$}.
            \For{$i = 1,\ldots,r$}
                \State Draw a random anchor \smash{$\boldsymbol{\theta} \in \prod_{j=1, j\not=t}^{d_\mathrm{ex}}\mathcal{D}_j$}.
                \State Sample function values $g_{\{t\}}(\xi^{(j)}, \boldsymbol{\theta})$.
                \State Perform sparse reconstruction of $\ghat_{\{t\},l}(\boldsymbol{\theta})$ for all $l \in \mathcal{P}_{\{t\}}(\Gamma)$.
                \State Add the $s$ largest indices $l$ satisfying $|\ghat_{\{t\},l}(\boldsymbol{\theta})| \geq \delta_+$ to $I_{\{t\}}$.
            \EndFor
        \EndFor
        
        \Statex \textbf{(Step 2) Coupled component identification}
        \For{$t=2,\ldots,J$}
            \State Set $\tilde r\coloneqq r$ if $t<J$ and $\tilde r\coloneqq 1$ otherwise.
            \State Set $K\coloneqq(I_{\{1,\ldots,t-1\}} \times I_{\{t\}})\cap\mathcal{P}_{\{1,\ldots,t\}}(\Gamma)$.
            \State Draw random sampling nodes \smash{$(\b \xi^{(j)})_{j=1}^M \subset \prod_{j=1}^t\mathcal{D}_j$}.
            \State Set $I_{\{1,\ldots,t\}}\coloneqq\emptyset$.
            \For{$i=1,\ldots,\tilde r$}
                \State Draw a random anchor \smash{$\boldsymbol{\theta} \in \prod_{j=t+1}^{d_\mathrm{ex}}\mathcal{D}_j$}.
                \State Sample function values $g_{\{1,\ldots,t\}}(\b \xi^{(j)}, \boldsymbol{\theta})$ (or  $g_{\{1,\ldots,t\}}(\b \xi^{(j)})$ if $t=J$).
                \State Perform sparse reconstruction of $\ghat_{\{1,\ldots,t\},\b l}(\boldsymbol{\theta})$ for all ${\b l \in K}$.
                \State Add the $s$ largest indices $\b l$ where $|\ghat_{\{1,\ldots,t\},\b k}(\boldsymbol{\theta})| \geq \delta_+$ to $I_{\{1,\ldots,t\}}$.
            \EndFor
        \EndFor
        
        \Statex \textbf{(Step 3)}
        \State Set $I_{\mathrm{est}}\coloneqq I_{\{1,\ldots,J\}}$ and $\ghat_{\b l} \coloneqq \ghat_{\{1,\ldots,J\},\b l}$ for all $\b l\in I$.

        \Ensure Detected index set $I_{\mathrm{est}}\subset\Gamma$, approximated coefficients $(\ghat_{\b l})_{\b l\in I} \in\mathbb{C}^{|I|}$.
    \end{algorithmic}
\end{algorithm}

\subsection{Orthogonal Matching Pursuit} \label{subsec:OMP}

Now, we briefly discuss OMP, which is used to reconstruct the relevant projected coefficients in lines 7 and 19 of Algorithm \ref{alg:main}.
This algorithm is a greedy approach developed for sparse approximation and compressed sensing.
For a broad overview, we refer to \cite{PaReKr93,TrGi07,KuRa08}.
Assume, we are in the $t$-th dimension-incremental step of Algorithm \ref{alg:main} and consider the candidate set $K$ as given in \eqref{eq:candidate_set} and sampling values $g_{\{1,\ldots,t\}}(\b  \xi) \coloneqq (g_{\{1,\ldots,t\}}(\b \xi^{(j)}))_{j=1}^M \in \R^M$.
Starting from an empty active set, OMP iteratively selects at each step the product of basis functions~$\phi_{\{1,\ldots,t\}, \b l}$ for which the corresponding column of the basis matrix
\begin{align}\label{eq:OMP_basis_mat}
\b\Phi \coloneqq (\phi_{\{1,\ldots,t\},\b l}(\b \xi^{(j)})_{j=1,\ldots,M,\, \b l \in K} \in \C^{M \times \vert K \vert}
\end{align}
is most strongly correlated with the current residual.
The active set is then updated with the index~$\b l$, and the basis coefficients $\ghat_{\{1,\ldots,t\},\b l}(\boldsymbol{\theta})$ restricted to this set are determined via a least squares solve.
This procedure is repeated until the prescribed sparsity~$s$ is reached or a suitable stopping criterion, e.g., an absolute or relative error tolerance, is satisfied.

Under suitable assumptions on the sampling points~\smash{$(\b\xi^{(j)})_{j=1}^M$}, the resulting Algorithm \ref{alg:OMP} is guaranteed to recover the correct support of the basis coefficients~$\ghat_{\{1,\ldots,t\},\b l}(\boldsymbol{\theta})$, $\b l \in K,$ with high probability provided that the number of samples scales like~\mbox{$M \gtrsim s \log |K|$}.
The complexity of a single iteration can be bounded by~\mbox{$\mathcal{O}(M |K| + M s + s^2)$}, cf.~\cite{StCh12}, which is dominated by the cost of the matrix-vector multiplication~($M |K|$).
Performing~$s$ iterations leads to the overall computational complexity~\mbox{$\mathcal{O}(s M |K|)$}.
Our implementation is based on the methodology presented in~\cite{MoNePoSoUl25}, with a Python implementation being publicly available\footnote{\url{https://github.com/Zeppo1994/SparseRecovery}}.

\begin{remark}
While we only discuss OMP here, other sparse recovery methods can be also employed.
Examples include the CoSaMP \cite{NeeTr09} and square-root LASSO \cite{BeCheWa11}.
Our Python implementation \cite{gitNABOPB} includes these approaches as additional choices.
In our numerical experiments, no method consistently outperformed the others such that we have chosen OMP for simplicity.
\end{remark}

\begin{algorithm}[t]
\caption{OMP}\label{alg:OMP}
\textbf{Input}: Sampled BOPB matrix \(\b\Phi \), sample vector $g_{\{1,\ldots,t\}}(\boldsymbol \xi) \in \R^M$, number of iterations $s$
\begin{algorithmic}[1]
\State \smash{$S^{(0)} = \emptyset$}
\ForAll{$k=0,\ldots,s-1$}
\State \smash{$ {\b l}^{(k+1)} = \operatorname{argmax}_{\b l \in K}\{\lvert({\b\Phi }^\ast ({\b\Phi }{\b z}^{(k)} - g_{\{1,\ldots,t\}}(\boldsymbol \xi))_{\b l}\rvert\}$}
\State \smash{\(S^{(k+1)} = S^{(k)} \cup \{ {\b l}^{(k+1)}\}\)}
\State \smash{${\b z}^{(k+1)} = \operatorname{argmin}_{\b z\in \mathbb{C}^{\lvert J\rvert}}\{\lVert {\b\Phi }\b z - g_{\{1,\ldots,t\}}(\boldsymbol \xi) \rVert_{\ell_2} \quad \text{s.t. } \operatorname{supp}(\b z) \subset S^{(k+1)}\}$}\phantomsection\label{OMP:ApproxUpdate}
\EndFor
\State \textbf{return} \smash{$\b z^{(s)}$}
\end{algorithmic}
\end{algorithm}

\subsection{Improving Sample Efficiency (OMP+)}
\label{subsec:virtual_points}

In Algorithm~\ref{alg:main}, each evaluation of $u$ at a node $(\b x^{(j)}, \tau^{(j)}, \b a^{(j)}) \in \mathcal D$ requires a separate numerical solve of \eqref{eq:general} if the parameter instances $\b a^{(j)}$ differ.
As described so far, Algorithm~\ref{alg:main} incorporates random samples $\{\b\xi^{(j)}\}_{j=1}^M$ conditional to a random anchor, and thus in principle only a single value \smash{$u(\b x^{(j)}, \tau^{(j)})$} from each PDE solve, even though the solution is typically available at many grid points $(\b x, \tau)$.
Given the substantial cost of PDE solves, this represents a significant inefficiency.
Fortunately, in contrast to the cubature-based approach in \cite{PoTa25}, OMP imposes no structural constraints on the sampling nodes $\{\b\xi^{(j)}\}_{j=1}^M \in \R^{t}$.
Thus, we can systematically incorporate multiple samples from a computed solution, as detailed below.

Consider a $t$-th dimension-incremental step with $t > d+1$ in Algorithm \ref{alg:main}.
Then, we have \smash{$\b\xi^{(j)} = (\b x^{(j)}, \tau^{(j)}, \b a^{(j)}_1)$} with $\b a^{(j)}_1 \in \mathbb{R}^{t - d - 1}$ and the anchor $\boldsymbol{\theta} = \b a_2 \in \mathbb{R}^{n + d + 1 - t}$.
Instead of evaluating the PDE solution $u(\cdot,\cdot,\b a^{(j)})$ for \smash{$\b a^{(j)} = (\b a^{(j)}_1, \b a_2)$}, $j = 1,\ldots,M$, only at the sampled $(\b x^{(j)}, \tau^{(j)})$, we additionally evaluate them at $b^{d+1}$ randomly chosen points $(\b x,\tau) \in \mathbb{R}^{d+1}$, which of course differ for each $\b a^{(j)}$.
We recommend $b=3$ as default.
The number of evaluations per solve must be carefully chosen to balance the exploration of the parameter and space-time domains.
A smaller amount of space-time evaluations requires a higher number of costly PDE solves to achieve a given target sample size.
Conversely, too many evaluations per parameter result in insufficient coverage of the parameter domain (the extreme case being only a single parameter).
Moreover, since the complexity of OMP scales linearly with $M$, adding more space-time evaluations is only recommended when they contain non-redundant information.
In our experiments, we refer to this modified version of \texttt{OMP} as \texttt{OMP+}.

\begin{remark}
\label{rem:add_mod}
In all one-dimensional steps of Algorithm \ref{alg:main}, where a coordinate of $(\b x, \tau)$ is the active variable, as well as in all dimension-incremental steps with $t \leq d+1$, the random anchor $\boldsymbol{\theta}$ fully covers the parameters $\b a$.
In these cases, a single PDE solve suffices to generate the samples.
While this applies to any reconstruction method, we only incorporate it into \texttt{OMP+}.
The latter already distinguishes between spatial, temporal, and parameter coordinates, making the modification straightforward to implement.
In contrast, \texttt{OMP} follows the original non-intrusive dimension-incremental framework, which does not have access to structural information.
\end{remark}

\section{Numerics}
\label{sec:numerics}

We benchmark Algorithm~\ref{alg:main} for several PDEs, in both its
standard version \texttt{OMP} and the modified variant \texttt{OMP+}
from Section~\ref{subsec:virtual_points}.
For $\b a \in [-1,1]^n$, we measure the approximation accuracy in terms of the absolute $\ell_2$-error
\begin{align}\label{eq:abs_err}
\text{err}_{\text{abs}}(\b a) &\coloneqq \norm{S_{I_{\mathrm{est}}}^\A u(\b x,\tau,\b a) - u(\b x,\tau,\b a)}_{\ell_2} = \biggl(\sum_{j=1}^{G}\abs{S_{I_{\mathrm{est}}}^\A u(\b x^{(j)},\tau^{(j)},\b a) - u(\b x^{(j)},\tau^{(j)},\b a)}^2\biggr)^{\frac12},
\end{align}
which coincides with the quality criterion in Algorithm \ref{alg:main}, and the (scale-invariant) relative $\ell_2$-error
\begin{align} \label{eq:rel_err}
\text{err}_{\text{rel}}(\b a) &\coloneqq \frac{\norm{S_{I_{\mathrm{est}}}^\A u(\b x,\tau,\b a) - u(\b x,\tau,\b a)}_{\ell_2}}{\norm{u(\b x,\tau,\b a)}_{\ell_2}} \notag\\
&= \frac{\left(\sum_{j=1}^{G}\abs{S_{I_{\mathrm{est}}}^\A u(\b x^{(j)},\tau^{(j)},\b a) - u(\b x^{(j)},\tau^{(j)},\b a)}^2\right)^{\frac12}}{\left(\sum_{j=1}^{G}\abs{u(\b x^{(j)},\tau^{(j)},\b a)}^2\right)^{\frac12}}.
\end{align}
In \eqref{eq:abs_err} and \eqref{eq:rel_err}, the $(\b x^{(j)},\tau^{(j)})$, $j=1,\ldots,G$, are equidistant grid points in the space-time domain $\Omega \times [0,T]$ (or spatial domain $\Omega$ for stationary problems).
We statistically compare these quantities for several randomly drawn coefficients $\b a$ in terms of range (min and max), median, first quartile, and third quartile.

Throughout, we employ the tensorized Chebyshev polynomial basis, which is given by
\begin{align}
T_{\b k}(\b z) \coloneqq \prod_{j=1}^{d_\mathrm{ex}} T_{k_j}(z_j) \quad\text{with}\quad  T_{k_j}(z_j) \coloneqq \begin{cases}
    1 & k_j=0 \\
    \sqrt2 \cos(k_j \arccos(z_j)) & k_j \not=0
\end{cases},
\end{align}
where $d_\mathrm{ex} = d+1+n$ denotes the dimension of the space-time-parameter domain as in Section~\ref{subsec:sparse_approximation}. 
This choice is motivated by the fact that several classes of parameter-dependent PDEs admit sparse or approximately sparse representations in terms of Chebyshev polynomials \cite{PoTa25}.
If not stated otherwise, Algorithm \ref{alg:main} uses the following parameters:

\begin{itemize}
\item the search space $\Gamma$ is the (non-negative) full grid $[0,N]^{d+1+n}$ in $d+1+n$ dimensions (or $d+n$ for stationary examples) with extension $N$;
\item the detection threshold $\delta_+ = 10^{-12}$;
\item$r=5$ detection iterations for \texttt{OMP} and $r=2$ for \texttt{OMP+}, see also Remark \ref{rem:detection_iterations};
\item an individual sparsity level $s$ for each experiment.
\end{itemize}
See \cite{KaPoTa23} for more detailed information on these parameters and how they affect the behavior of the support detection.
For the one-dimensional detection steps, we use $s$ samples.
This suffices because typically $N \ll s$, which ensures that the number of candidates in these steps (namely $N+1$) remains much smaller than $s$.
For the coupled detection steps, the standard \texttt{OMP} approach utilizes an oversampling factor of $c=5$, leading to $5s$ samples and PDE solves.
In contrast, the \texttt{OMP+} approach utilizes $3^{d+1}$ (or $3^d$ for stationary problems) additional space-time evaluations as discussed in Section~\ref{subsec:virtual_points}, reducing the requirement to only $s$ PDE solves.
In both cases, the resulting sparse recovery problems \eqref{eq:SparseRecProb} are solved using the OMP Algorithm~\ref{alg:OMP} detailed in Section~\ref{subsec:OMP}.

We compare our methods against the dimension-incremental method using Chebyshev multiple rank-1 lattices (\texttt{cMR1L}) \cite[Sec.\ 3]{PoTa25} that inspired our work.
Further, we consider Tensorized Fourier Neural Operators (\texttt{TFNO}), implemented in PyTorch using the NeuralOperator library \cite{gitNeuralOperator}.
Given training samples \smash{$(u_0^{(j)}, u^{(j)})$}, the underlying Fourier Neural Operator (FNO) framework learns the solution operator $G \colon u_0 \mapsto u$ from data by minimizing the supervised Sobolev loss
\begin{align} \label{eq:TFNO_min}
\min_{\theta} \;\frac{1}{M} \sum_{j=1}^M 
\|G_\theta(u_0^{(j)}) - u^{(j)}\|_{H^1}^2.
\end{align}
The neural operator $G_\theta$ is realized with a pointwise lifting layer, followed by a composition of multiple mappings (frequently called Fourier layers) of the form
\begin{align}
v_{k+1}(\b x)
=
\sigma \left(
W_k v_k(\b x)
+
\mathcal{F}^{-1}
\big(
R_k \cdot \mathcal{F}(v_k)
\big)(\b x)
\right),
\end{align}
where $R_k$ are learned Fourier multipliers on a truncated set of modes and $W_k$ are pointwise linear mappings, and a final pointwise projection layer.
For more details on FNOs and their tensorized variants, we refer to~\cite{LiEtAl21,kossaifi2024multigrid}.
All \texttt{TFNO} models use $4$ Fourier layers with $64$ hidden channels and are trained using the AdamW optimizer with initial learning rate $10^{-3}$ and weight decay $10^{-4}$.
Training is performed for 1000 epochs with a batch size of 64 using random parameter samples $\b a_i$ with corresponding solutions $u_{\b a_i}$ evaluated at resolution $100$.
Further, we incorporate a cosine annealing schedule that decays the learning rate to $10^{-8}$ and is stepped after each epoch.

All experiments are performed in Python\textsuperscript{TM} and can be found together with the algorithm in \cite{gitNABOPB}.
Unless stated otherwise, all computations are carried out on two AMD EPYC 9534 64-Core processors.
To accelerate the sampling process, parallelization with up to $96$ workers is employed.
The \texttt{TFNO} models are trained separately on a NVIDIA GeForce RTX 4070 Ti Super GPU.
Individual runtimes are discussed in the corresponding sections.

\subsection{Heat equation}\label{subsec:heat}

Our first numerical example is the heat equation with homogeneous boundary conditions on the domain $\Omega = (0,1)$ with final time $T=1$, namely
\begin{equation} \label{eq:heat_1d}
\begin{alignedat}{3}
\partial_\tau u - \alpha^2 \partial_{xx} u &= 0, & \qquad x \in (0,1),\, \tau&\in (0,1) \\
u(x,0) &= f(x), & \qquad x &\in (0,1)\\
u(0,\tau) = u(1,\tau) &= 0, & \qquad \tau &\in (0,1).
\end{alignedat}
\end{equation}
For the experiments, we choose $\alpha = 0.25$.
As in \cite[Sec.~3.5]{PoTa25}, we use initial conditions given by a truncated sine expansion, namely
\begin{align}\label{eq:heat_initial}
u(x,0) = f(x) = \sum_{\ell=1}^9 a_\ell \sin (\ell \pi x), & & x \in [0,1]
\end{align}
with coefficients $a_\ell \in [-1,1], \ell \in \N$, which leads to the analytical solution
\begin{align}\label{eq:heat_exact}
u(x,\tau) = \sum_{\ell=1}^9 a_\ell \sin \left(\ell \pi x\right) \exp \left(-\ell^2\pi^2\alpha^2 \tau \right), & & x \in [0,1], \tau\in [0,1].
\end{align}

Clearly, the truncation together with the condition $a_\ell \in [-1,1]$ restricts the set of admissible initial conditions $f$.
Functions far outside this class require a modified parameterization, for instance by increasing the number of retained modes the range of the coefficients.
As explained in Remark \ref{rem:transform}, we shift both PDE variables $x$ and $t$ via the invertible transformation $\mathcal{T}\colon [-1,1]\to[0,1]$ with $\mathcal{T}(z) =
\frac12(z+1)$ to be able to use the Chebyshev basis.

\begin{remark}
Choosing multivariate Chebyshev polynomials as the underlying BOPB can be theoretically justified for several of our examples.
For the heat equation, the explicit solution formula \eqref{eq:heat_exact} shows that each input $a_\ell$ enters the solution only linearly, and only through a single summand.
This linear dependence is reproduced exactly by the degree-one Chebyshev
polynomial $T_1(a_\ell)=\sqrt2\,a_\ell$ (with $T_0\equiv1$ in the remaining
coordinates), so that \eqref{eq:heat_exact} may be rewritten as
\begin{align}
u(x,\tau,\b a)
&=
\sum_{\ell=1}^{9}
\frac{T_1(a_\ell)}{\sqrt2}
\sin(\ell\pi x)
\exp\!\bigl(-\ell^2\pi^2\alpha^2\tau\bigr)
\\
&=
\sum_{\ell=1}^{9}
\frac{T_{\b e_\ell}(\b a)}{\sqrt2}
\sin(\ell\pi x)
\exp\!\bigl(-\ell^2\pi^2\alpha^2\tau\bigr),
\end{align}
where $\b e_\ell=(0,\ldots,0,1,0,\ldots,0)^\top\in\N_0^9$ is the $\ell$-th standard unit vector.
Consequently, the relevant parameter dependence is concentrated on only the
nine multi-indices $\b e_1,\ldots,\b e_9$.
One therefore expects a suitable index set $I_{\mathrm{est}}$ to consist of
indices of the form $\b k=(\ast,\ast,\b e_\ell)^\top$, $\ell=1,\ldots,9$,
where the first two entries (the indices in the $x$- and $\tau$-directions)
remain arbitrary.
Hence, the corresponding support is not only sparse but also highly structured.
\end{remark}

For the sampling, we solve the PDE \eqref{eq:heat_1d} using SciPy's \texttt{solve\_ivp} with the implicit Radau IIA Runge--Kutta method of order five, a relative tolerance of $10^{-8}$, and an absolute tolerance of $10^{-10}$.
In Algorithm \ref{alg:main}, we use sparsity $s=1000$ and extension $N=64$ for the search space $\Gamma$, i.e., $\Gamma=[0,64]^{11}$ and $|\Gamma|\approx 8.75 \cdot 10^{19}$.

The \texttt{cMR1L} comparison from \cite[Sec.~3.5]{PoTa25} uses the same parameters.
As a learning-based baseline, we train a TFNO with $529\,233$ parameters on $10\,000$ and $40\,000$ training samples, respectively, resulting in $10\,100$ and $40\,500$ total PDE solves including validation data.
We denote these approaches as \texttt{TFNO}$_{10k}$ and \texttt{TFNO}$_{40k}$, respectively.


\begin{figure}[hp]
\centering
\begin{tikzpicture}
\begin{axis}[
	boxplot/draw direction=x,
	xmode=log,
	xmin=0.000009,
	xmax=0.005,
	ytick={1,2,3,4,5,6,7},
	yticklabels={\texttt{cMR1L},\texttt{OMP},\texttt{OMP+},\texttt{TFNO}$_{40k}$,\texttt{TFNO}$_{10k}$},
	xlabel={$\text{err}_{\text{rel}}(\b a)$},
	width=0.95\textwidth,
    height=0.30\textwidth,
	]
	\pgfplotstableread[col sep=comma]{heat_boxerr.csv}\datatable;
	\addplot+[
  		boxplot prepared from table={
    			table=\datatable,
    			row=0,
    			lower whisker=lw,
    			upper whisker=uw,
    			lower quartile=lq,
    			upper quartile=uq,
    			median=med
  		}, my boxplot style,
    ] coordinates {};
    \addplot+[
  		boxplot prepared from table={
    			table=\datatable,
    			row=1,
    			lower whisker=lw,
    			upper whisker=uw,
    			lower quartile=lq,
    			upper quartile=uq,
    			median=med
  		}, my boxplot style,
    ] coordinates {};
    \addplot+[
  		boxplot prepared from table={
    			table=\datatable,
    			row=2,
    			lower whisker=lw,
    			upper whisker=uw,
    			lower quartile=lq,
    			upper quartile=uq,
    			median=med
  		}, my boxplot style,
    ] coordinates {};
    \addplot+[
  		boxplot prepared from table={
    			table=\datatable,
    			row=3,
    			lower whisker=lw,
    			upper whisker=uw,
    			lower quartile=lq,
    			upper quartile=uq,
    			median=med
  		}, my boxplot style,
    ] coordinates {};
    \addplot+[
  		boxplot prepared from table={
    			table=\datatable,
    			row=4,
    			lower whisker=lw,
    			upper whisker=uw,
    			lower quartile=lq,
    			upper quartile=uq,
    			median=med
  		}, my boxplot style,
    ] coordinates {};
\end{axis}
\end{tikzpicture}
\caption{Relative approximation error $\text{err}_{\text{rel}}(\b a)$ for $10000$ randomly drawn $\b a$ for the heat equation.
Box-and-whisker plots show median, first and third quartiles, and maximum and minimum errors.}
\label{fig:heat_err}

\vspace{1cm}

\centering
\begin{tikzpicture}
\begin{axis}[
    ybar stacked,
    bar width=50pt,
    width=0.95\textwidth,
    height=0.40\textwidth,
    ylabel={runtime in minutes},
    xtick={1,2,3,4,5},
    xticklabels={\texttt{cMR1L},\texttt{OMP},\texttt{OMP+},\texttt{TFNO}$_{40k}$,\texttt{TFNO}$_{10k}$},
    enlarge x limits=0.22,
    ymajorgrids=true,
    grid style=dashed,
    ymin=0,
    ymax=370,
    ytick={0,60,120,180,240},
    yticklabels={0,60,120,180,240},
    legend style={
        at={(0.5,1.02)},
        anchor=south,
        legend columns=-1,
        /tikz/every odd column/.append style={column sep=3pt},
        /tikz/every even column/.append style={column sep=20pt},
        draw=none
    },
    extra y ticks={280,320,360},
    extra y tick labels={355,700,1420},
    extra y tick style={grid=none},
    legend image code/.code={
    \draw[#1, draw=black] (0cm,-0.1cm) rectangle (0.6cm,0.3cm);
    },
]

\draw[dashed,black,thick] (0,240) -- (6,240);

\addplot[
    fill=white,
    draw=black,
    pattern=north east lines,
    pattern color=black,
] coordinates {
    (1, 319.60)   
    (2,   193.51)
    (3,   28.88)
    (4, 33.92)
    (5,  8.91)
};

\addplot[
    fill=white,
    draw=black,
    pattern=dots,
    pattern color=black,
] coordinates {
    (1, 0.40)  
    (2,   7.57)
    (3, 29.01)
    (4, 326.08)
    (5, 271.09)
};

\coordinate (labcmr) at (axis cs:1,145);
\coordinate (labomp) at (axis cs:2,100);
\coordinate (labompp) at (axis cs:3,60);
\coordinate (labtfno40) at (axis cs:4,50);
\coordinate (labtfno10) at (axis cs:5,30);
\coordinate (gputfno40) at (axis cs:4.15,300);
\coordinate (gputfno10) at (axis cs:5.15,260);

\legend{runtime needed for solves, remaining runtime}
\end{axis}
\node[fill=white, draw=none, font=\footnotesize, align=center]
    at (labcmr) {$702\,054$\\solves};

\node[fill=white, draw=none, font=\footnotesize, align=center]
    at (labomp) {$285\,000$\\solves};

\node[fill=white, draw=none, font=\footnotesize, align=center, anchor=south]
    at (labompp) {$35\,006$\\solves};

\node[fill=white, draw=none, font=\footnotesize, align=center, anchor=south]
    at (labtfno40) {$40\,500$\\solves};

\node[fill=white, draw=none, font=\footnotesize, align=center, anchor=south]
    at (labtfno10) {$10\,100$\\solves};

\node[
    fill=white,
    draw=black!50,
    rounded corners=1pt,
    inner sep=1pt,
    font=\small
] at (gputfno40) {GPU};

\node[
    fill=white,
    draw=black!50,
    rounded corners=1pt,
    inner sep=1pt,
    font=\small
] at (gputfno10) {GPU};

\draw[thick] ([xshift=-2pt]current axis.west |- {rel axis cs:0,0.68})
      -- ([xshift=2pt]current axis.west |- {rel axis cs:0,0.70});
\draw[thick] ([xshift=-2pt]current axis.west |- {rel axis cs:0,0.71})
      -- ([xshift=2pt]current axis.west |- {rel axis cs:0,0.73});
\draw[thick] ([xshift=-2pt]current axis.east |- {rel axis cs:0,0.68})
      -- ([xshift=2pt]current axis.east |- {rel axis cs:0,0.70});
\draw[thick] ([xshift=-2pt]current axis.east |- {rel axis cs:0,0.71})
      -- ([xshift=2pt]current axis.east |- {rel axis cs:0,0.73});
\end{tikzpicture}
\caption{Runtime and sample complexity for the heat equation.
The upper axis is compressed for visualization.
Each bar is split into PDE solve time (lower) and remaining runtime (upper). 
All computations were performed on the same CPU, except for \texttt{TFNO} training, which used a GPU.}
\label{fig:heat_runtime}

\vspace{1cm}

\centering
\begin{tikzpicture}
\begin{axis}[
	boxplot/draw direction=x,
	xmode=log,
	xmin=0.000004,
	xmax=1,
	ytick={1,2,3,4,5},
	yticklabels={$250$,$500$,$750$,$1000$,$1500$},
	xlabel={$\text{err}_{\text{rel}}(\b a)$},
    ylabel={sparsity $s$},
	width=0.95\textwidth,
    height=0.30\textwidth,
	]
	\pgfplotstableread[col sep=comma]{heat_boxerr.csv}\datatable;
    \addplot+[
  		boxplot prepared from table={
    			table=\datatable,
    			row=5,
    			lower whisker=lw,
    			upper whisker=uw,
    			lower quartile=lq,
    			upper quartile=uq,
    			median=med
  		}, my boxplot style,
    ] coordinates {};
    \addplot+[
  		boxplot prepared from table={
    			table=\datatable,
    			row=6,
    			lower whisker=lw,
    			upper whisker=uw,
    			lower quartile=lq,
    			upper quartile=uq,
    			median=med
  		}, my boxplot style,
    ] coordinates {};
    \addplot+[
  		boxplot prepared from table={
    			table=\datatable,
    			row=7,
    			lower whisker=lw,
    			upper whisker=uw,
    			lower quartile=lq,
    			upper quartile=uq,
    			median=med
  		}, my boxplot style,
    ] coordinates {};
    \addplot+[
  		boxplot prepared from table={
    			table=\datatable,
    			row=2,
    			lower whisker=lw,
    			upper whisker=uw,
    			lower quartile=lq,
    			upper quartile=uq,
    			median=med
  		}, my boxplot style,
    ] coordinates {};
    \addplot+[
  		boxplot prepared from table={
    			table=\datatable,
    			row=8,
    			lower whisker=lw,
    			upper whisker=uw,
    			lower quartile=lq,
    			upper quartile=uq,
    			median=med
  		}, my boxplot style,
    ] coordinates {};
\end{axis}
\end{tikzpicture}
\caption{Relative approximation error $\text{err}_{\text{rel}}(\b a)$ for $10000$ randomly drawn $\b a$ for the heat equation using \texttt{OMP+} with varying sparsity $s$.
Box-and-whisker plots show median, first and third quartiles, and maximum and minimum errors.}
\label{fig:heat_OMP+_err}
\end{figure}

Figure~\ref{fig:heat_err} illustrates the relative approximation error for all methods.
\texttt{OMP+} achieves a median error of approximately $4\cdot 10^{-5}$, comparable to \texttt{cMR1L}.
This improves over \texttt{OMP}, which attains a median error of roughly $2\cdot 10^{-4}$.
The baselines \texttt{TFNO$_{40k}$} and \texttt{TFNO$_{10k}$} yield slightly larger median errors of about $3\cdot 10^{-4}$ and $5\cdot 10^{-4}$, respectively.

Figure~\ref{fig:heat_runtime} shows the number of solves, the corresponding computation times, and the remaining runtime excluding PDE solves.
\texttt{OMP} required slightly more than $3$ hours, with the dominant portion spent on solving the PDEs for sample generation.
Similarly, \texttt{cMR1L} required roughly $700$ hours due to the significantly larger number of required samples.
The PDE-specific \texttt{OMP+} version was the fastest overall, with a total runtime of about $1$ hour, roughly half of which was spent on solving the PDEs.
The \texttt{TFNO$_{10k}$} and \texttt{TFNO$_{40k}$} required approximately $6$ and $23$ hours of training time, respectively.

The effect of the sparsity parameter $s$ on the approximation quality is illustrated in Figure~\ref{fig:heat_OMP+_err} for \texttt{OMP+} with extension $N=64$.
The relative error $\text{err}_{\text{rel}}$ decreases as $s$ increases, while the spread of errors remains within roughly one order of magnitude in each setting.

\subsection{Burgers' equation}

Our second example is the non-linear 1D Burgers' equation with
homogeneous boundary conditions, namely
\begin{equation} \label{eq:burger_1d}
\begin{alignedat}{3}
\partial_\tau u + u \partial_x u &= \nu \partial_{xx} u, & \qquad x \in (0,1),\, \tau &\in (0,1) \\
u(x,0) &= f(x), & \qquad x &\in (0,1)\\
u(0,\tau) = u(L,\tau) &= 0, & \qquad \tau &\in (0,1),
\end{alignedat}
\end{equation}
which we solve for the viscosity parameter $\nu = 0.05$. 
In contrast to the previous example, we are only interested in the solution $u(\cdot,1)$ at the final time.
Again, given $\b a \in [-1,1]^9$, we consider initial conditions $f$ given by the truncated sine expansion
\begin{align}\label{eq:burger_initial_approx}
f(x) \approx \sum_{\ell=1}^9 a_\ell \sin (\ell \pi x), & & x \in [0,1].
\end{align}
As before, we transform the spatial variable via $\mathcal{T} x = \frac12(x+1)$ and use \texttt{solve\_ivp} with the Radau method, this time with a relative tolerance of $10^{-6}$ and an absolute tolerance of ${10^{-8}}$.
In Algorithm \ref{alg:main}, we also use sparsity $s=1000$ and extension $N=64$.
Since the target function space is one-dimensional, \texttt{OMP+} would
use only $3$ space-time evaluations per PDE solve in the recommended default configuration.
To improve the efficiency, we increase this to $9$, giving $10$ evaluations per solve, consistent with the previous example.
Regarding the results in \cite{PoTa25}, already sparsity $s=500$ and extension $N=32$ lead to runtimes on the order of several days, making larger configurations infeasible.
We highlight this mismatch in $s$ by denoting the comparison as \texttt{cMR1L}$^\ast$.
The \texttt{TFNO} baseline with $81\,049$ parameters was trained on $40\,000$ samples, resulting in $40\,500$ total PDE solves including validation data.


\begin{figure}[hp]
\centering
\begin{subfigure}{0.95\textwidth}
\centering
\begin{tikzpicture}
\begin{axis}[
	boxplot/draw direction=x,
	xmode=log,
	xmin=0.00005,
	xmax=10,
	ytick={1,2,3,4},
	yticklabels={\texttt{cMR1L}$^\ast$,\texttt{OMP},\texttt{OMP+},\texttt{TFNO}},
	xlabel={$\text{err}_{\text{rel}}(\b a)$},
	width=0.95\textwidth,
    height=0.30\textwidth,
	]
	\pgfplotstableread[col sep=comma]{burger_boxerr.csv}\datatable;
	\addplot+[
  		boxplot prepared from table={
    			table=\datatable,
    			row=0,
    			lower whisker=lw,
    			upper whisker=uw,
    			lower quartile=lq,
    			upper quartile=uq,
    			median=med
  		}, my boxplot style,
    ] coordinates {};
    \addplot+[
  		boxplot prepared from table={
    			table=\datatable,
    			row=1,
    			lower whisker=lw,
    			upper whisker=uw,
    			lower quartile=lq,
    			upper quartile=uq,
    			median=med
  		}, my boxplot style,
    ] coordinates {};
    \addplot+[
  		boxplot prepared from table={
    			table=\datatable,
    			row=2,
    			lower whisker=lw,
    			upper whisker=uw,
    			lower quartile=lq,
    			upper quartile=uq,
    			median=med
  		}, my boxplot style,
    ] coordinates {};
    \addplot+[
  		boxplot prepared from table={
    			table=\datatable,
    			row=3,
    			lower whisker=lw,
    			upper whisker=uw,
    			lower quartile=lq,
    			upper quartile=uq,
    			median=med
  		}, my boxplot style,
    ] coordinates {};
    \end{axis}
\end{tikzpicture}
\caption{The relative approximation error $\text{err}_{\text{rel}}(\b a)$.}
\end{subfigure}

\vspace{0.4cm}

\begin{subfigure}{0.95\textwidth}
\centering
\begin{tikzpicture}
\begin{axis}[
	boxplot/draw direction=x,
	xmode=log,
	xmin=0.00009,
	xmax=1,
	ytick={1,2,3,4},
	yticklabels={\texttt{cMR1L}$^\ast$,\texttt{OMP},\texttt{OMP+},\texttt{TFNO}},
	xlabel={$\text{err}_{\text{abs}}(\b a)$},
	width=0.95\textwidth,
    height=0.30\textwidth,
	]
	\pgfplotstableread[col sep=comma]{burger_boxerr.csv}\datatable;
	\addplot+[
  		boxplot prepared from table={
    			table=\datatable,
    			row=6,
    			lower whisker=lw,
    			upper whisker=uw,
    			lower quartile=lq,
    			upper quartile=uq,
    			median=med
  		}, my boxplot style,
    ] coordinates {};
    \addplot+[
  		boxplot prepared from table={
    			table=\datatable,
    			row=5,
    			lower whisker=lw,
    			upper whisker=uw,
    			lower quartile=lq,
    			upper quartile=uq,
    			median=med
  		}, my boxplot style,
    ] coordinates {};
    \addplot+[
  		boxplot prepared from table={
    			table=\datatable,
    			row=4,
    			lower whisker=lw,
    			upper whisker=uw,
    			lower quartile=lq,
    			upper quartile=uq,
    			median=med
  		}, my boxplot style,
    ] coordinates {};
    \addplot+[
  		boxplot prepared from table={
    			table=\datatable,
    			row=7,
    			lower whisker=lw,
    			upper whisker=uw,
    			lower quartile=lq,
    			upper quartile=uq,
    			median=med
  		}, my boxplot style,
    ] coordinates {};
    \end{axis}
\end{tikzpicture}
\caption{The absolute approximation error $\text{err}_{\text{abs}}(\b a)$.}
\end{subfigure}
\caption{The relative and absolute approximation error $\text{err}(\b a)$ for $10000$ randomly drawn $\b a$ for the Burgers' equation example. The box-and-whisker plots show the median, the first and the third quartile as well as the maximal and minimal error observed.
Note that \texttt{cMR1L}$^\ast$ uses a smaller sparsity $s$ and extension $N$.}
\label{fig:burger_err}

\vspace{1cm}

\centering
\begin{tikzpicture}
\begin{axis}[
    ybar stacked,
    bar width=50pt,
    width=0.95\textwidth,
    height=0.40\textwidth,
    ylabel={runtime in minutes},
    xtick={1,2,3,4},
    xticklabels={\texttt{cMR1L}$^\ast$,\texttt{OMP},\texttt{OMP+},\texttt{TFNO}},
    enlarge x limits=0.18,
    ymajorgrids=true,
    grid style=dashed,
    ymin=0,
    ymax=380,
    ytick={0,60,120,180,240,300},
    yticklabels={0,60,120,180,240,300},
    legend style={
        at={(0.5,1.02)},
        anchor=south,
        legend columns=-1,
        /tikz/every odd column/.append style={column sep=3pt},
        /tikz/every even column/.append style={column sep=20pt},
        draw=none
    },
    extra y ticks={360},
    extra y tick labels={6858},
    extra y tick style={grid=none},
    legend image code/.code={
        \draw[#1, draw=black] (0cm,-0.1cm) rectangle (0.6cm,0.3cm);
    },
]

\draw[dashed,black,thick] (0,300) -- (5,300);

\addplot[
    fill=white,
    draw=black,
    pattern=north east lines,
    pattern color=black,
] coordinates {
    (1, 360.0)      
    (2, 175.50)
    (3, 50.17)
    (4, 55.19)
};

\addplot[
    fill=white,
    draw=black,
    pattern=dots,
    pattern color=black,
] coordinates {
    (1, 1.10)      
    (2,   112.49)
    (3, 12.93)
    (4, 207.81)
};

\coordinate (labcmr) at (axis cs:1,180);
\coordinate (labomp) at (axis cs:2,100);
\coordinate (labompp) at (axis cs:3,70);
\coordinate (labtfno) at (axis cs:4,70);
\coordinate (gputfno) at (axis cs:4.11,243);

\legend{runtime needed for solves, remaining runtime}
\end{axis}

\draw[thick] ([xshift=-2pt]current axis.west |- {rel axis cs:0,0.84})
      -- ([xshift=2pt]current axis.west |- {rel axis cs:0,0.86});
\draw[thick] ([xshift=-2pt]current axis.west |- {rel axis cs:0,0.87})
      -- ([xshift=2pt]current axis.west |- {rel axis cs:0,0.89});
\draw[thick] ([xshift=-2pt]current axis.east |- {rel axis cs:0,0.84})
      -- ([xshift=2pt]current axis.east |- {rel axis cs:0,0.86});
\draw[thick] ([xshift=-2pt]current axis.east |- {rel axis cs:0,0.87})
      -- ([xshift=2pt]current axis.east |- {rel axis cs:0,0.89});

\node[fill=white, draw=none, font=\footnotesize, align=center]
    at (labcmr) {$5\,706\,215$\\solves};

\node[fill=white, draw=none, font=\footnotesize, align=center]
    at (labomp) {$255\,000$\\solves};

\node[fill=white, draw=none, font=\footnotesize, align=center, anchor=south]
    at (labompp) {$35\,002$\\solves};

\node[fill=white, draw=none, font=\footnotesize, align=center, anchor=south]
    at (labtfno) {$40\,500$\\solves};

\node[
    fill=white,
    draw=black!50,
    rounded corners=1pt,
    inner sep=1pt,
    font=\small
] at (gputfno) {GPU};

\end{tikzpicture}
\caption{Runtime and sample complexity for the Burgers' equation.
The upper axis is compressed for visualization.
Each bar is split into PDE solve time (lower) and remaining runtime (upper).
All computations were performed on the same CPU, except for \texttt{TFNO} training, which used a GPU.}
\label{fig:burger_runtime}
\end{figure}

Figure \ref{fig:burger_err} shows the relative and absolute approximation errors for all approaches.
The \texttt{TFNO} median error is almost two orders of magnitude smaller than those of the remaining approaches, which range between $10^{-2}$ and $10^{-1}$.
The relative error $\text{err}_{\text{rel}}(\b a)$ exhibits large upper whiskers for all approaches, which are absent for the absolute error $\text{err}_{\text{abs}}(\b a)$.

Figure \ref{fig:burger_runtime} shows the computation time and number of PDE solves for each approach.
\texttt{cMR1L}$^\ast$ required significantly more solves than the rest, resulting in substantially longer runtime.
\texttt{OMP} and \texttt{TFNO} both finished in less than $5$ hours, while \texttt{OMP+} needed only about $1$ hour.

\begin{figure}[t]
\centering
\begin{tikzpicture}
\begin{axis}[
	boxplot/draw direction=x,
	xmode=log,
	xmin=0.004,
	xmax=1,
	ytick={1,2,3,4,5,6},
	yticklabels={{$3$ SE,\,$c=1$},{$9$ SE,\,$c=1$},{$15$ SE,\,$c=1$},{$3$ SE,\,$c=3$},{$9$ SE,\,$c=3$},{$15$ SE,\,$c=3$}},
	xlabel={$\text{err}_{\text{abs}}(\b a)$},
	width=0.95\textwidth,
    height=0.30\textwidth,
	]
	\pgfplotstableread[col sep=comma]{burger_boxerr.csv}\datatable;
	\addplot+[
  		boxplot prepared from table={
    			table=\datatable,
    			row=8,
    			lower whisker=lw,
    			upper whisker=uw,
    			lower quartile=lq,
    			upper quartile=uq,
    			median=med
  		}, my boxplot style,
    ] coordinates {};
    \addplot+[
  		boxplot prepared from table={
    			table=\datatable,
    			row=4,
    			lower whisker=lw,
    			upper whisker=uw,
    			lower quartile=lq,
    			upper quartile=uq,
    			median=med
  		}, my boxplot style,
    ] coordinates {};
    \addplot+[
  		boxplot prepared from table={
    			table=\datatable,
    			row=11,
    			lower whisker=lw,
    			upper whisker=uw,
    			lower quartile=lq,
    			upper quartile=uq,
    			median=med
  		}, my boxplot style,
    ] coordinates {};
    \addplot+[
  		boxplot prepared from table={
    			table=\datatable,
    			row=9,
    			lower whisker=lw,
    			upper whisker=uw,
    			lower quartile=lq,
    			upper quartile=uq,
    			median=med
  		}, my boxplot style,
    ] coordinates {};
    \addplot+[
  		boxplot prepared from table={
    			table=\datatable,
    			row=10,
    			lower whisker=lw,
    			upper whisker=uw,
    			lower quartile=lq,
    			upper quartile=uq,
    			median=med
  		}, my boxplot style,
    ] coordinates {};
    \addplot+[
  		boxplot prepared from table={
    			table=\datatable,
    			row=12,
    			lower whisker=lw,
    			upper whisker=uw,
    			lower quartile=lq,
    			upper quartile=uq,
    			median=med
  		}, my boxplot style,
    ] coordinates {};
    \end{axis}
\end{tikzpicture}
\caption{Absolute approximation error $\text{err}(\b a)$ for $10000$ randomly drawn $\b a$ for the Burgers' equation using \texttt{OMP+} with sparsity $s=1000$, varying numbers of additional spatial evaluations (SE) and oversampling factors $c$.
Box-and-whisker plots show median, first and third quartiles, and maximum and minimum errors.}
\label{fig:burger_err_OMP+}
\end{figure}
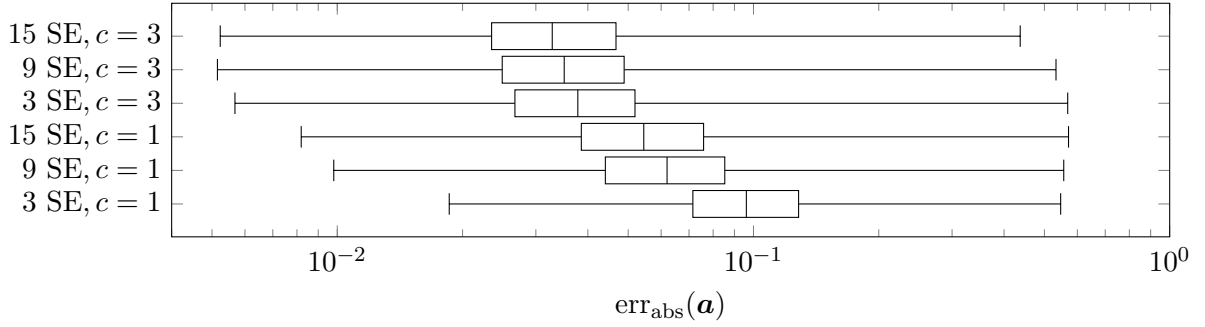

Lastly, Figure \ref{fig:burger_err_OMP+} depicts the absolute approximation error $\text{err}_{\text{abs}}(\b a)$ for \texttt{OMP+} with varying numbers of additional spatial evaluations per PDE solve ($3$, $9$, and $15$) and oversampling factors ($c \in \{1, 3\}$).
The error magnitude remains similar across all six variations, with accuracy improving slightly for more spatial evaluations or higher oversampling factor $c$.

\subsection{Parametric diffusion equation}

As third example, we consider the 2D diffusion equation on $\Omega = [0,1]^2$ with homogeneous Dirichlet boundary conditions and an affine random diffusion coefficient in some domain $\Omega_{\b y}$, see \cite{BaCoMi16,BaCoVo17} namely
\begin{equation}\label{eq:rand_diff_eq}
\begin{aligned}
-\nabla \cdot (a(\b x,\b y)\nabla u(\b x,\b y)) &= f(\b x), \qquad & \b x \in \Omega,\ \b y \in \Omega_{\b y},\\
u(\b x,\b y) &= 0, \qquad & \b x \in \partial\Omega,\ \b y \in \Omega_{\b y}.
\end{aligned}
\end{equation}
Here, the differential operator $\nabla$ acts with respect to the spatial variable $\b x$.
As in \cite[Sec.~11]{EiGiSchwZa14} and \cite[Sec.~4.3]{KaPoTa22}, we use the affine coefficient
\begin{align}
a(\b x,\b y) = 1 + \sum_{j=1}^{20} y_j \psi_j(\b x), \qquad \b x \in \Omega,\ \b y \in [-1,1]^{20},
\end{align}
with 
$\psi_j(\b x) = c j^{-\mu} \cos(2\pi m_1(j)x_1)\cos(2\pi m_2(j)x_2).
$
The parameter are chosen as $c = 0.9/\zeta(2)$ and $\mu = 2$, together with $k(j) \coloneqq \lfloor -1/2 + \sqrt{1/4 + 2j} \rfloor$, $m_1(j) \coloneqq j-\frac{k(j) (k(j)+1)}{2}$ and $m_2(j) \coloneqq k(j)-m_1(j)$.
Moreover, we set $f \equiv 1$.
As before, we transform $\b x$ componentwise by the affine map $\mathcal T(z)=\frac12(z+1)$.
In total, this yields a $22$-dimensional approximation problem.

As in \cite{KaPoTa22}, we solve the PDE \eqref{eq:rand_diff_eq} numerically using FEniCS on a mesh with $5000$ finite elements.
In Algorithm \ref{alg:main}, we choose $s=1000$ and $N=64$.
Moreover, we restrict the candidate set $\Gamma$ by imposing $\max_{\b k \in \Gamma}\|\b k\|_0 = 7$.
Such a restriction reflects the common observation that, in many high-dimensional parametric PDEs, higher-order interactions contribute considerably less than lower-order ones.
This choice is adopted from \cite[Sec.~3.4]{PoTa25} to ensure a comparable experimental setup.
The \texttt{cMR1L} comparison uses the same parameters.
The \texttt{TFNO} with $529\,233$ parameters is trained on $85\,000$ samples, resulting in a total of $86\,000$ PDE solves including validation data.

Figure~\ref{fig:diff_pde_err} compares the relative approximation error $\mathrm{err}(\b y)$ for all four approaches.
The \texttt{TFNO} achieves a median error of roughly $2.5\cdot10^{-4}$, almost an order of magnitude below the other approaches, whose median errors are around $10^{-3}$.
On the other hand, the \texttt{TFNO}'s upper whisker spans almost a full order of magnitude, compared to roughly half an order of magnitude for the other approaches.

Figure~\ref{fig:diff_pde_runtime} shows the number of PDE solves and corresponding runtimes for all approaches.
As in the previous example, \texttt{cMR1L} required significantly more solves and runtime than \texttt{OMP} and \texttt{OMP+}, which finished in roughly $8$ and $2$ hours, respectively.
Notably, the \texttt{TFNO} training time is comparable to the total runtime of \texttt{cMR1L}.

\begin{figure}[t]
\centering
\begin{tikzpicture}
\begin{axis}[
	boxplot/draw direction=x,
	xmode=log,
	xmin=0.0001,
	xmax=0.01,
	ytick={1,2,3,4},
	yticklabels={\texttt{cMR1L},\texttt{OMP},\texttt{OMP+},\texttt{TFNO}},
	xlabel={$\text{err}_{\text{rel}}(\b a)$},
	width=0.95\textwidth,
    height=0.30\textwidth,
	]
	\pgfplotstableread[col sep=comma]{diff_pde_boxerr.csv}\datatable;
	\addplot+[
  		boxplot prepared from table={
    			table=\datatable,
    			row=0,
    			lower whisker=lw,
    			upper whisker=uw,
    			lower quartile=lq,
    			upper quartile=uq,
    			median=med
  		}, my boxplot style,
    ] coordinates {};
    \addplot+[
  		boxplot prepared from table={
    			table=\datatable,
    			row=1,
    			lower whisker=lw,
    			upper whisker=uw,
    			lower quartile=lq,
    			upper quartile=uq,
    			median=med
  		}, my boxplot style,
    ] coordinates {};
    \addplot+[
  		boxplot prepared from table={
    			table=\datatable,
    			row=2,
    			lower whisker=lw,
    			upper whisker=uw,
    			lower quartile=lq,
    			upper quartile=uq,
    			median=med
  		}, my boxplot style,
    ] coordinates {};
    \addplot+[
  		boxplot prepared from table={
    			table=\datatable,
    			row=3,
    			lower whisker=lw,
    			upper whisker=uw,
    			lower quartile=lq,
    			upper quartile=uq,
    			median=med
  		}, my boxplot style,
    ] coordinates {};
    \end{axis}
\end{tikzpicture}
\caption{Relative approximation error $\text{err}_{\text{rel}}(\b y)$ for $10000$ randomly drawn $\b y$ for the parametric diffusion equation.
Box-and-whisker plots show median, first and third quartiles, and maximum and minimum errors.}
\label{fig:diff_pde_err}

\vspace{1cm}

\centering
\begin{tikzpicture}
\begin{axis}[
    ybar stacked,
    bar width=50pt,
    width=0.95\textwidth,
    height=0.40\textwidth,
    ylabel={runtime in minutes},
    xtick={1,2,3,4},
    xticklabels={\texttt{cMR1L},\texttt{OMP},\texttt{OMP+},\texttt{TFNO}},
    enlarge x limits=0.18,
    ymajorgrids=true,
    grid style=dashed,
    ymin=0,
    ymax=840,
    ytick={0,120,240,360,480,600},
    yticklabels={0,120,240,360,480,600},
    legend style={
        at={(0.5,1.02)},
        anchor=south,
        legend columns=-1,
        /tikz/every odd column/.append style={column sep=3pt},
        /tikz/every even column/.append style={column sep=20pt},
        draw=none
    },
    extra y ticks={720,780},
    extra y tick labels={3466,3570},
    extra y tick style={grid=none},
    legend image code/.code={
        \draw[#1, draw=black] (0cm,-0.1cm) rectangle (0.6cm,0.3cm);
    },
]

\draw[dashed,black,thick] (0,600) -- (5,600);

\addplot[
    fill=white,
    draw=black,
    pattern=north east lines,
    pattern color=black,
] coordinates {
    (1, 717.99)  
    (2, 441.11)
    (3, 84.13)
    (4, 103.24)
};

\addplot[
    fill=white,
    draw=black,
    pattern=dots,
    pattern color=black,
] coordinates {
    (1, 2.01)      
    (2, 70.22)
    (3, 39.03)
    (4, 677) 
};

\coordinate (labcmr) at (axis cs:1,300);
\coordinate (labomp) at (axis cs:2,200);
\coordinate (labompp) at (axis cs:3,140);
\coordinate (labtfno) at (axis cs:4,140);
\coordinate (gputfno) at (axis cs:4.11,730);

\legend{runtime needed for solves, remaining runtime}
\end{axis}

\draw[thick] ([xshift=-2pt]current axis.west |- {rel axis cs:0,0.76})
      -- ([xshift=2pt]current axis.west |- {rel axis cs:0,0.78});
\draw[thick] ([xshift=-2pt]current axis.west |- {rel axis cs:0,0.79})
      -- ([xshift=2pt]current axis.west |- {rel axis cs:0,0.81});
\draw[thick] ([xshift=-2pt]current axis.east |- {rel axis cs:0,0.76})
      -- ([xshift=2pt]current axis.east |- {rel axis cs:0,0.78});
\draw[thick] ([xshift=-2pt]current axis.east |- {rel axis cs:0,0.79})
      -- ([xshift=2pt]current axis.east |- {rel axis cs:0,0.81});

\node[fill=white, draw=none, font=\footnotesize, align=center]
    at (labcmr) {$3\,683\,946$\\solves};

\node[fill=white, draw=none, font=\footnotesize, align=center]
    at (labomp) {$615\,000$\\solves};

\node[fill=white, draw=none, font=\footnotesize, align=center, anchor=south]
    at (labompp) {$79\,006$\\solves};

\node[fill=white, draw=none, font=\footnotesize, align=center, anchor=south]
    at (labtfno) {$86\,000$\\solves};

\node[
    fill=white,
    draw=black!50,
    rounded corners=1pt,
    inner sep=1pt,
    font=\small
] at (gputfno) {GPU};

\end{tikzpicture}
\caption{Runtime and sample complexity for the parametric diffusion equation. The upper axis is compressed for visualization.
Each bar is split into PDE solve time (lower) and remaining runtime (upper). 
All computations were performed on the same CPU, except for \texttt{TFNO training}, which used a GPU.}
\label{fig:diff_pde_runtime}
\end{figure}

\subsection{Discussion}\label{subsec:discussion}

Now, we discuss the numerical results and the observed trade-offs between approximation quality and computational effort.

\paragraph{Sparse approximation}

A key observation is that the approximation accuracy of \texttt{OMP} and \texttt{OMP+} depends primarily on the sparsity of the solution $u$ in the chosen basis.
For the heat equation, $u$ is highly sparse in the multivariate Chebyshev basis due to the linear dependence on the
parameters $a_\ell$ in \eqref{eq:heat_exact}.
Consequently, the index set $I_{\mathrm{est}}$ detected by Algorithm \ref{alg:main} already performs well for small sparsity $s = |I_{\mathrm{est}}|$, and increasing the sparsity $s$ further improves accuracy, as seen in Figure \ref{fig:heat_OMP+_err}.
The increase saturates once the diameter of the search space $\Gamma$ (which is restricted by $N$) is exhausted.

The Burgers' equation, by contrast, showed the opposite behavior: $u$ is not well represented by sparse Chebyshev expansions, as discussed in \cite[Sec.~3.6]{PoTa25}.
This results in much lower approximation accuracy at the same sparsity $s$.
Moreover, increasing the amount of training data improves accuracy only marginally, as seen in Figure \ref{fig:burger_err_OMP+}, suggesting that the observed limitations stem from the basis rather than insufficient data for the sparse recovery step.

These results highlight that the choice of basis critically affects the approximation quality.
When the detected index set $I_{\mathrm{est}}$ exhibits little sparsity structure, experimenting with alternative bases is advisable.
Indeed, \texttt{OMP} and \texttt{OMP+} are applicable for arbitrary product bases, provided that an efficient routine for the matrix-vector multiplication with the basis matrix \eqref{eq:OMP_basis_mat} is available.

\paragraph{Samples and PDE solves}

Applying sparse recovery methods such as OMP to the dimension-incremental framework yields a substantial advantage in sample complexity over cubature approaches based on rank-1 lattices.
While the rank-1 lattice methods in \cite{PoTa25} enable efficient reconstruction via FFT, this benefit is outweighed by the large number of required PDE solves, which scales as $s^2$ in the worst case.
In contrast, the number of samples required by \texttt{OMP} scales only linearly in $s$.
Furthermore, the additional samples per PDE solve in  \texttt{OMP+} (which are not usable for cubature-based methods) yield a substantial additional reduction in computational cost (by up to $80\%$) without sacrificing approximation quality, clearly visible across all experiments.
As a consequence, our methods close the gap in terms of sample complexity
to learning-based approaches such as TFNOs, though the dimension-incremental
framework does require samples with a particular structure.

Moreover, we want to emphasize that the \texttt{TFNO} models serve purely as baselines, following the generic guidelines of the neural operator library \cite{gitNeuralOperator} without further modification.
The comparison should therefore not be interpreted as a rigorous benchmark, but rather as an indication of the approximation quality and computational complexity achievable by our methods relative to off-the-shelf neural operator configurations.

\paragraph{Interpretability of the approximation}

An important advantage of the proposed  approach is its interpretability.
The detected index set $I_{\mathrm{est}}$ reveals which input variables are most influential, which parameter interactions are relevant, and to what extent higher-order interactions matter.
The resulting approximation is thus not merely a black-box predictor, but also yields qualitative insight into the underlying PDE.
In contrast, such structural information is far less accessible in neural operator approaches such as TFNOs, where the learned representation is distributed across a large number of trainable parameters.

\bibliographystyle{abbrv}

\end{document}